\documentclass[12pt, oneside,reqno]{amsart}
\usepackage[dvipsnames]{xcolor}
\usepackage{amsfonts,amsmath,mathabx,fullpage,amssymb,amsthm,tikz,mathtools,lmodern,paralist,mathrsfs,hyperref,stmaryrd,enumitem,tikz-cd,ytableau,tcolorbox,cases, soul, accents}
\usepackage[normalem]{ulem}
\hypersetup{colorlinks, citecolor=red, filecolor=black, linkcolor=blue, urlcolor=blue}
\usetikzlibrary{calc,3d,patterns}

\usepackage{color}
\usepackage[color]{xy}
\colorlet{cyan}[rgb]{cyan} 
\ytableausetup{smalltableaux, aligntableaux=center}

\newcommand{\excise}[1]{}

\DeclareMathOperator{\Comp}{\mathsf{Comp}}
\DeclareMathOperator{\WComp}{\mathsf{WComp}}
\DeclareMathOperator{\codim}{codim}

\DeclareMathOperator{\link}{link}
\DeclareMathOperator{\des}{des}

\DeclareMathOperator{\Symm}{Sym}

\DeclareMathOperator{\shuffle}{\mathsf{Shuffle}}
\DeclareMathOperator{\initial}{ini}

\DeclareMathOperator{\anti}{\mathsf{s}}
\DeclareMathOperator{\Scr}{\mathsf{Scr}}

\newcommand{\defterm}[1]{\textbf{\boldmath{#1}\unboldmath}}

\newcommand{\kk}{\Bbbk}

\newcommand{\Rr}{\mathbb{R}}

\newcommand{\Zz}{\mathbb{Z}}

\newcommand{\BB}{\mathcal{B}}
\newcommand{\CC}{\mathcal{C}}

\newcommand{\II}{\mathcal{I}}

\newcommand{\MM}{\mathcal{M}}
\newcommand{\NN}{\mathcal{N}}

\newcommand{\QQ}{\mathcal{Q}}

\newcommand{\UU}{\mathcal{U}}

\newcommand{\CCC}{\mathsf{C}}

\newcommand{\HS}{\mathbf{H}_\QQ}

\newcommand{\hatzero}{\hat{\mathbf{0}}}

\newcommand{\inj}{\hookrightarrow}

\newcommand{\xx}{\mathbf{x}}
\newcommand{\yy}{\mathbf{y}}
\newcommand{\zz}{\mathbf{z}}

\newcommand{\pp}{\mathfrak{p}}
\newcommand{\qq}{\mathfrak{q}}

\newcommand{\0}{\emptyset}
\newcommand{\isom}{\cong}
\newcommand{\sm}{\setminus}
\newcommand{\ov}{\overline}
\newcommand{\Sym}{\mathfrak{S}}
\newcommand{\x}{\times}
\newcommand{\compn}{\vDash}
\newcommand{\wcompn}{\vDash_w}

\newcommand{\un}{\mathbf{1}_\kk}

\newcommand{\setcomp}[1]{\uppercase{#1}}

\newcommand{\refineseq}{\trianglerighteq}
\newcommand{\refinedbyeq}{\trianglelefteq}

	\newcommand{\bH}{\mathbf{H}}
	
\newcommand{\bl}{\pmb{\ell}}		\newcommand{\bL}{\mathbf{L}}
\newcommand{\GP}{\mathbf{GP}}
\newcommand{\OGP}{\mathbf{OGP}}

\newcommand{\Mat}{\mathbf{Mat}}
\newcommand{\OMat}{\mathbf{OMat}}

\numberwithin{equation}{section}

\usepackage{thmtools}

\declaretheorem[style=plain,numberwithin=section]{theorem}

\declaretheorem[style=plain,sibling=theorem]{corollary}

\declaretheorem[style=plain,sibling=theorem]{proposition}

\declaretheorem[style=definition,sibling=theorem]{definition}
\declaretheorem[style=definition,qed=$\blacktriangleleft$,sibling=theorem]{example}

\newcommand{\dllin}{\ar@{-}[dl]}
\newcommand{\drlin}{\ar@{-}[dr]}
\newcommand{\dlin}{\ar@{-}[d]}

\begin{document}
\title{Hopf monoids of ordered simplicial complexes}

\author{Federico Castillo}
\address{Department of Mathematics \\ University of Kansas}
\email{fcastillo@ku.edu}
\author{Jeremy L.\ Martin}
\address{Department of Mathematics \\ University of Kansas}
\email{jlmartin@ku.edu}
\author{Jos\'e A.\ Samper}
\address{Max-Planck Institute for Mathematics in the Sciences}
\email{samper@mis.mpg.de}
\date{\today}
\maketitle
\begin{abstract}We study ordered matroids and generalized permutohedra from a Hopf theoretic point of view. Our main object is a Hopf monoid in the vector species of extended generalized permutahedra equipped with an order of the coordinates; this monoid extends the Hopf monoid of generalized permutahedra studied by Aguiar and Ardila.  Our formula for the antipode is cancellation-free and multiplicity-free, and is supported only on terms that are compatible with the local geometry of the polyhedron.  Our result is part of a larger program to understand orderings on ground sets of simplicial complexes (for instance, on shifted and matroid independence complexes).  In this vein, we show that shifted simplicial complexes and broken circuit complexes generate Hopf monoids that are expected to exhibit similar behavior.
\end{abstract}
\section{Introduction}

Shifted complexes and matroid independence complexes are well-known families of simplicial complexes that share many nice properties, such as shellability and Laplacian integrality.  Both kinds of complexes can be studied through orderings of their vertex sets.  The very definition of a shifted complex requires a fixed vertex ordering, which then controls its structure strongly.  Dually, matroids can be informally characterized as complexes for which the choice of ordering does not matter for constructions such as Kruskal's algorithm, the activity formula for the Tutte polynomial, or the homotopy type of the broken-circuit complex.

The goal of this work is to investigate orderings of shifted complexes, matroids and related objects from a Hopf-theoretic point of view.  Matroids can be considered as special cases of generalized permutahedra, and Aguiar and Ardila's study \cite{AA} of the Hopf monoid $\GP$ of generalized permutahedra is a significant starting point for our work.  In this context, a linear order on the ground set of a matroid generalizes to an order of coordinates in the ambient space of a polyhedron.



Accordingly, we define an \defterm{ordered generalized permutahedron} to be a pair $(w,\pp)$ (or a tensor $w\otimes\pp$), where $w$ is a linear order on a finite set $I$ and $\pp\subset\Rr^I$ is a generalized permutahedron.  The corresponding Hopf monoid $\OGP$ is the Hadamard product $\bL^*\x\GP$, where $\bL^*$ is the Hopf monoid of linear orderings equipped with shuffle product and $\GP$ is the monoid of generalized permutahedra studied by Aguiar and Ardila.  We use $\bL^*$ rather than $\bL$ (linear orderings with concatenation product) because $\bL^*$ and $\OGP$, hence $\bL^*\x\GP$, are commutative, but $\bL$ is not.

Our first main result (Theorem~\ref{OGP-antipode}) is a formula for the antipode $\anti(w\otimes\pp)$ in $\OGP$.  The formula is cancellation-free (all terms involve distinct basis elements) and multiplicity-free (all coefficients are $0$ or $\pm 1$).  As in the Aguiar--Ardila formula for the antipode in $\GP$, the coefficients can be interpreted topologically, as Euler characteristics of simplicial complexes arising from normal cones of faces of~$\pp$.  In the course of this calculation, we are led to study a family of simplicial complexes that we call \defterm{Scrope complexes}, which have a simple combinatorial description and may be of independent interest.  Furthermore, the antipode formula is \textit{local}, in the sense that the support of $\anti(w\otimes\pp)$ is restricted to tensors $u\otimes\qq$ such that $\qq$ is a face of $\pp$ containing the vertex whose normal cone contains the braid chamber associated with $w$.

The original motivation of studying orderings of simplicial complexes leads us to define a \defterm{Hopf class}, which is a class of ordered simplicial complexes closed under suitable ordered analogues of join, restriction, and contraction (the simplicial operations required to define product and coproduct).  The definition collapses in the unordered setting (underlining the need to work with ordered complexes), and every Hopf class indeed gives rise to a Hopf monoid (Theorem~\ref{thm:quasiHopf}).  Somewhat surprisingly, there is also a ``universal'' Hopf class $\UU$ that contains all Hopf classes as subclasses, and contains all matroid and pure shifted complexes, as well as broken-circuit complexes (Theorem~\ref{universal-Hopf}).  These results suggest that ordered simplicial complexes should be understood from a geometric and Hopf-theoretic point of view. 

\section{Background and notation} \label{sec:background}

We begin by setting up definitions and notation for the objects we will need, including preposets, polytopes, generalized permutahedra, and Hopf monoids.  Our presentation owes a great deal to \cite{AA} and \cite{PRW}, although our notation and terminology differs from theirs in some cases.
Generalized permutahedra were introduced by Postnikov \cite{Beyond}.

\subsection{Simplicial complexes and matroids} \label{sec:scs-and-matroids}
A \defterm{simplicial complex} $\Sigma$ on a finite set $I$ is a (possibly empty) subset of $2^I$ closed under inclusion. The elements of $I$ are \defterm{vertices}, the elements of $\Sigma$ are \defterm{faces} and the maximal (under inclusion) faces are \defterm{facets}.  A collection of faces $\sigma_1,\dots,\sigma_r\subseteq I$ \defterm{generate} the complex $\langle\sigma_1,\dots,\sigma_r\rangle$, namely the union of their power sets.  
The complex $\Sigma$ is \defterm{pure} if all facets have the same size.  The \defterm{reduced Euler characteristic} of $\Sigma$ is $\tilde\chi(\Sigma)=\sum_{\sigma\in\Sigma} (-1)^{|\sigma|-1}$.  This coincides with the reduced Euler characteristic of the standard topological realization $|\Sigma|$.  An \defterm{ordered complex} on a finite set $I$ is a pair $(w,\Sigma)$, where $w$ is a total ordering on $I$ and $\Sigma$ is a simplicial complex on $I$. The following notions are fundamental in both the ordered and unordered settings, so we introduce them together.

Let $S\subseteq I$. The \defterm{restriction} $\Sigma|_S$ is the complex with vertex set $S$ and faces $\{S\cap\sigma:\sigma\in\Sigma\}$.  (This complex may also be referred to as the \defterm{deletion} of $I\sm S$, or the subcomplex \defterm{induced by $S$}.) The \defterm{link} of $\sigma\in\Sigma$ is $\link_\Sigma(\sigma)=\{\tau\in\Sigma:\ \sigma\cap\tau=\emptyset,\ \sigma\cup\tau\in\Sigma\}$; this is a simplicial complex on $I\sm S$. In the ordered setting, a contraction of a subset $S$ of $I$ is the ordered complex $(\Sigma\slash S, w|_{I\backslash S})$ where $\Sigma\slash S$ is the link of the smallest facet of $\Sigma|_S$ in the lexicographic order induced by $w$.

If $\Sigma_1,\Sigma_2$ are simplicial complexes on disjoint vertex sets $I_1,I_2$, then their \defterm{join} $\Sigma_1*\Sigma_2$ is $\Sigma_1*\Sigma_2=\{\sigma_1\cup\sigma_2:\ \sigma_1\in \Sigma_1,\ \sigma_2\in \Sigma_2\}$.

A pure simplicial complex $\Sigma$ on $I$ is a \defterm{matroid independence complex}, or simply a \defterm{matroid}, if $\Sigma|_S$ is pure for every $S\subseteq I$.
In standard matroid theory terminology, $I$ is usually called the \defterm{ground set} of the matroid, faces of $\Sigma$ are called \defterm{independent sets}, and facets are called \defterm{bases}.  The \defterm{direct sum} of two matroids is just their join as simplicial complexes.  The \defterm{contraction} $\Sigma/I$ can be defined as the link of any facet of $\Sigma|_I$.  Note that all restrictions and contractions of a matroid complex are themselves matroid complexes.

\subsection{Set compositions, preposets, and the braid fan} \label{sec:set-comps-and-gps}

Let $S$ be a finite set.  A \defterm{(set) composition} of $S$ is an ordered list of nonempty, pairwise-disjoint subsets $A_1,\dots,A_k\subseteq S$ (\defterm{blocks}) whose union is $S$.  A \defterm{weak (set) composition} is defined similarly except that
the blocks are allowed to be empty.  The number of blocks is $k=|A|$.  The symbols $\Comp(S)$ and $\WComp(S)$ denote the sets of compositions and weak compositions of $S$, and we abbreviate $\Comp(n)=\Comp([n])$ and $\WComp(n)=\WComp([n])$.  We also write $A\compn S$ or $A\wcompn S$ to indicate respectively $A\in\Comp(S)$ or $A\in\WComp(S)$.  In both cases we typically write $A=A_1|\cdots|A_k$.  In this notation, the vertical bars are called \defterm{separators}.  Note that the order of elements within each block is not significant.  A set of compositions of~$S$ is called an \defterm{album}.

The set $\Comp(n)$ is partially ordered by refinement: $A\refineseq B$ means that every block of $B$ is of the form $A_i\cup A_{i+1}\cup\cdots\cup A_{j-1}\cup A_j$. 
In this case, the equivalence relation $\equiv_B$ imposes an equivalence relation on the blocks of $A$.  
The refinement ordering is ranked and has a unique minimal element, namely the composition $\hatzero$ with one block. Every linear order $w$ gives rise to a set composition $\setcomp{w}$ with $n$ singleton blocks, namely
$\setcomp{w}=w(1)\,|\,w(2)\,|\,\cdots\,|\,w(n)$.

\begin{definition} \label{defn:descent-comp}
Let $w,u$ be linear orders on $I$.  The \defterm{$u$-descent composition of $w$} is the set composition $D_{w,u}$ of $\setcomp{w}:=w(1)|\cdots|w(|I|)$ for which $w(i)\equiv w(i+1)$ if and only if $(i,i+1)$ is a descent of $u^{-1}w$ (equivalently, if $w(i)$ occurs before $w(i+1)$ in $u$).
\end{definition}

For example, let $w=\mathsf{aebfcdhg}$, $u=\mathsf{bdahfgce}$ be linear orders of $I=\{\mathsf{a,b,c,d,e,f,g,h}\}$.  Then $u^{-1}w=38\cdot157\cdot246$ (with the descents marked).  Then $$D_{w,u}=\mathsf{ae|bfc|dhg=ae|bcf|dgh}.$$

A \defterm{preposet} $Q$ on $S$ is given by  a relation $\preceq_Q$ on~$S$ that is reflexive ($x\preceq_Qx$ for all $x\in S$) and transitive (if $x\preceq_Qy$ and $y\preceq_Qz$, then $x\preceq_Qz$).  The notation $x\equiv_Qy$ means that both $x\preceq_Qy$ and $x\succeq_Qy$; this is evidently an equivalence relation, whose equivalence classes are called the \defterm{blocks} of $Q$.  An \defterm{antichain} in $Q$ is a subset $T\subseteq S$ such that $x\not\prec y$ for all $x,y\in T$.  

The preposet $Q$ gives rise to a poset $Q/\!\!\equiv_Q$ on its blocks.  If this poset is a chain, then $Q$ is a \defterm{preorder}.  A \defterm{linear extension} of a preposet $Q$ is a preorder $R$ with the same blocks as $Q$ and such that $x\preceq_Qy$ implies $x\preceq_R y$ for all $x,y$.  A preorder $R$ contains the same information as a set composition given by the blocks.  

If the underlying set $S$ of a preposet $Q$ is equipped with a total order $<$ (e.g., if $S\subseteq\Zz$), then we say that $Q$ is \defterm{natural} with respect to $<$ if $x<y$ whenever $x\prec_Qy$.  If $Q$ is natural, then the poset $Q/\!\!\equiv_Q$ inherits a natural labeling from~$Q$.

\begin{figure}[ht]
\begin{center}
\begin{tikzpicture}
\newcommand{\sh}{5}
\draw (0,0)--(0,1) (1,1)--(1,0)--(2,1)--(2,0)--(3,1);
\node[fill=white] at (0,1) {\sf4};		\node[fill=white] at (1,1) {\sf3};		\node[fill=white] at (2,1) {\sf6};		\node[fill=white] at (3,1) {\sf8};
\node[fill=white] at (0,0) {\sf25};		\node[fill=white] at (1,0) {\sf1};		\node[fill=white] at (2,0) {\sf7};
\node at (-1.5,.5) {$\boxed{Q}$};
\begin{scope}[shift={(7,0)}]
\draw (0,1)--(0,0)--(1,1)--(1,2);
\node[fill=white] at (0,1) {\sf2345};	\node[fill=white] at (1,1) {\sf67};
\node[fill=white] at (0,0) {\sf1};		\node[fill=white] at (1,2) {\sf8};
\node at (2.5,.5) {$\boxed{Q^\natural}$};
\end{scope}
\end{tikzpicture}
\end{center}
\caption{Naturalization of a preposet.\label{fig:nat}}
\end{figure}
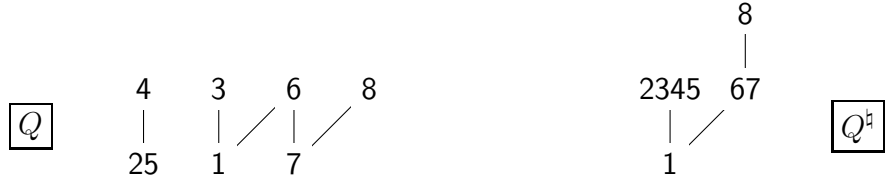

We will later need the notion of the \defterm{naturalization} of $Q$ with respect to a total order $<$, which is defined as the preposet $Q^\natural$ obtained from $Q$ by repeatedly (i) identifying every pair of blocks $B,C$ such that $B\prec_Q C$ and there exist $b\in B$ and $c\in C$ with $b>c$; and (ii) identifying every pair of blocks such that each one contains an element less than an element of the other.  Thus $Q^\natural$ is a natural preposet whose blocks are those of a natural composition.  An example of the construction is shown in Figure~\ref{fig:nat}.

The \defterm{closure} of a preposet $Q$ is the album
\begin{equation}\label{eq:album-closure}
\CCC_Q= \{A\in\Comp(n):\ i\preceq_Qj ~\implies~ i\preceq_Aj\}.
\end{equation}

We now relate these definitions to the geometry of the \defterm{braid arrangement}.  The faces of the braid arrangement are relatively-open cones that partition $\Rr^n$; the set of all faces is called the \defterm{braid fan} $\BB_n$.
Every composition $A\compn[n]$ determines a 
relatively open face $\sigma_A\in\BB_n$ with $\dim\sigma_A=|A|$, namely
\[\sigma_A = \{(x_1,\dots,x_n)\in\Rr^n:\ x_i \mathrel{\substack{<\\=\\>}} x_j\text{ if and only if } i \mathrel{\substack{\prec_A\\=_A\\\succ_A}} j \text{ respectively } \},\]
and this correspondence is a bijection.  In fact $A\refineseq B$ if and only if $\overline{\sigma_A}\supseteq\sigma_B$ (where the bar denotes topological closure), so the correspondence may be viewed as an isomorphism of posets.  In particular, the maximal faces $\sigma_w\in\BB_n$ correspond to permutations $w\in\Sym_n$.
For each preposet $Q$, the album $\CCC_Q$ corresponds to the closed subfan
\[\CC_Q = \{\sigma_A\in\BB_n:\ A\in\CCC_Q\}\]
whose maximal faces correspond to the linear extensions of $Q$.  The closed subfans of $\BB_n$ that arise in this way are precisely those whose union is convex.  In addition, $\CC_Q = \overline{\bigcup_A \sigma_A}$,
where $A$ ranges over all linear extensions of $Q$.

\subsection{Generalized permutahedra} \label{sec:background-for-gp}

Let $\pp\subset\Rr^n$ be a polyhedron.  For each $\xx\in\Rr^n$, let $\lambda_\xx$ be the linear functional on $\Rr^n$ given by $\lambda_\xx(\yy)=\xx\cdot\yy$, and let $\pp_\xx$ be the face of $\pp$ maximized by $\lambda_\xx$.  The \defterm{normal cone} of a face $\qq\subset\pp$ is
\[N^\circ_\pp(\qq)=\{\xx\in\Rr^n:\ \pp_\xx=\qq\}.\]
This is a relatively open polyhedral cone of dimension $n-\dim\qq$.  The normal cones of faces comprise the \defterm{normal fan} $\NN(\pp)$. The polytope $\pp$ is a \defterm{generalized permutahedron} (or GP) if and only if its normal fan is a coarsening of the braid fan.  Note that in this case $\dim\pp\leq n-1$, because $\pp$ must be contained in some hyperplane orthogonal to the line spanned by $(1,\dots,1)$ (the smallest face in the braid fan).

Every set composition $A\compn n$ gives rise to a face $\pp_A\subseteq\pp$ defined by
\begin{equation} \label{maximizing-face}
\pp_A=\{\xx\in\Rr^n:\ \lambda(\xx)\geq\lambda(\yy) \ \ \forall \lambda\in\sigma_A,\ \yy\in\pp\}.
\end{equation}
If $A$ is a maximal set composition (i.e., with $n$ blocks), then the braid cone $\sigma_A$ has full dimension, hence is contained in a full-dimensional cone of $\NN(\pp)$, so $\pp_A$ is a vertex of $\pp$ (and all vertices arise in this way).
Moreover, for each face $\qq\subseteq\pp$, the album of compositions
\begin{equation} \label{normal-preposet}
\{A\compn n:\ \sigma_A\subseteq N^\circ_\pp(\qq)\} = \{A\compn n:\ \pp_A=\qq\}
\end{equation}
consists precisely of the set compositions coarsening some preposet $Q$ on $[n]$, the \defterm{normal preposet of $\qq$}.  Often we will work simultaneously with a face $\qq$ and its normal preposet~$Q$, which contain equivalent information.

The definition of a generalized permutahedron implies that normal cones of faces carry combinatorial structure.  Accordingly, we define the following fans and their corresponding albums:
\[\begin{array}{r@{}l@{}l l r@{}l@{}l}
\CC_\qq^\circ	&{}= \{\sigma_A:~ \sigma_A\subseteq N_\pp(\qq)\}
			&{}= \{\sigma_A:~ \pp_A=\qq\}, &\quad&
\CCC_Q^\circ	&{}= \{A:~ \sigma_A\in\CC_\qq^\circ\}
			&{}= \{A:~ \pp_A=\qq\},\\
\CC_\qq		&{}= \{\sigma_A:~ \sigma_A\subseteq\ov{N_\pp(\qq)}\}
			&{}= \{\sigma_A:~ \pp_A\supseteq\qq\}, &&
\CCC_Q		&{}= \{A:~ \sigma_A\in\CC_\qq\}
			&{}= \{A:~ \pp_A\supseteq\qq\},\\
\partial\CC_\qq	&{}= \CC_\qq\,\sm\,\CC^\circ_\qq
			&{}= \{\sigma_A:~ \pp_A\supsetneq\qq\}, &&
\partial\CCC_Q	&{}= \CCC_Q\,\sm\,\CCC^\circ_Q
			&{}= \{A:~ \pp_A\supsetneq\qq\}.
\end{array}\]

An \defterm{extended generalized permutahedron} (EGP) \cite[Defn.~4.2]{AA} is a polyhedron whose normal fan coarsens some \emph{convex subfan} of the braid fan $\BB_n\subset\Rr^I=\Rr^n$; in particular, it is not bounded if the normal fan is not complete.  Most of the above statements about generalized permutahedra can be carried over to this more general setting.

\section{Hopf monoids}

A Hopf monoid is a structure that can be thought of algebraically as a generalization of a group, or combinatorially as a framework for putting together and taking part labeled objects of a particular type (graphs, posets, matroids, etc.)    A comprehensive treatment of Hopf monoids can be found in~\cite{AgMa}, and an accessible ``user's guide'' in~\cite[\S2]{AA}; here we review the essentials.

Consider the category whose objects are finite sets and whose morphisms are bijections. A \defterm{(vector) species} is a functor from this category to the category of vector spaces over a field $\kk$. The \defterm{unit species} $\un$ is the species such that $\un[\emptyset] = \kk$ and $\un[I]= 0$ for all non-empty sets. 

A \defterm{Hopf monoid} is a species $\bH$ together with a collection $\mu$ of linear maps $\mu_{I,J}: \bH[I]\otimes \bH[J] \to \bH[I \sqcup J]$ (\defterm{products}), a morphism $\iota: \un \to \bH$ (the \defterm{unit}), a collection $\Delta$ of linear maps $\Delta_{I,J} : \bH[I\sqcup J] \to \bH[I]\otimes \bH[J]$ (\defterm{coproducts}), and a morphism $\epsilon:\bH\to \un$ (the \defterm{counit}), satisfying several technical conditions of which the most important are associativity of $\mu$, coassociativity of $\Delta$ (defined by reversing the arrows in the diagram for associativity), and compatibility between product and coproduct (briefly, $\mu$ is a comonoid morphism and $\Delta$ is a monoid morphism). A hopf monoid is \defterm{connected} iff $\bH[\emptyset]=\kk$.

The \defterm{antipode} $\anti$ in $\bH$ is a collection of maps $\anti^{\bH}_I:\bH[I]\to\bH[I]$ defined by a certain commutative diagram that generalizes group inversion (regarding a Hopf monoid as a generalization of a group).  For connected Hopf monoids it is given explicitly by the Takeuchi formula
\begin{equation} \label{Takeuchi}
\anti^{\bH}_I := \sum_{A\in\WComp(I)} (-1)^{|A|}\mu_{A}\circ \Delta_{A}.
\end{equation}
The Takeuchi formula is general, but involves substantial cancellation, so when studying a particularHopf monoid it is desirable to give a cancellation-free form

\subsection{The Hopf monoid \texorpdfstring{$\bL^*$}{bL*} of linear orders}

The Hopf monoid $\bL^*$ is defined as follows.  As a vector species, $\bL^*[I]$ is the $\kk$-vector space spanned by the set $\bl[I]$ of linear orders of $I$.  To define the product and coproduct on $\bL^*$, we first need some combinatorial preliminaries.  First, let $w^{(1)},\dots,w^{(q)}$ be a collection of linear orders on pairwise-disjoint sets $I_1,\dots,I_q$.
A \defterm{shuffle} of $w^{(1)},\dots,w^{(q)}$ is an ordering on $I_1\cup\cdots\cup I_q$ that restricts to $w_j$ on each $I_j$.  The set of all shuffles is denoted $\shuffle(w^{(1)},\dots,w^{(q)})$.  For example,
$\shuffle(12,3) = \{123,132,312\}$ and $\shuffle(12,34) = \{1234,1324,1342,3124,3142,3412\}$.
The shuffle operation is commutative and associative. The product on $\bL^*$ is defined using shuffles:
\[\mu_{I,J}(w,u)=\sum_{v\in\shuffle(w,u)}v\]

Second, let $w=(w(1),\dots,w(n))\in\bl[I]$.  An \defterm{initial segment} of $w$ is a linear order of the form $(w(1),\dots,w(k))$, where $0\leq k\leq n$.  The set of all initial segments of $w$ is denoted $\initial(w)$.  With this in hand, the coproduct on $\bL^*$ is defined by
\[\Delta_{I,J}(v) = \begin{cases}
v|_I\otimes v|_J & \text{ if } v|_I \in \initial(w),\\
0 & \text { otherwise.}
\end{cases}\]

\subsection{The Hopf monoid \texorpdfstring{$\GP$}{GP} of generalized permutahedra}

Let $\GP[I]$ be the $\kk$-vector space spanned by all generalized permutahedra in $\Rr^I$.  To make the vector species $\GP$ into a Hopf monoid, we define a product and coproduct by
\begin{equation}
\mu_{I,J}(\pp_1 \otimes \pp_2) = \pp_1\times \pp_2,\qquad  \Delta_{I,J}(\pp) = \pp|I \otimes \pp/I.
\end{equation}
where $\pp|I$ and $\pp/I$ are defined in \cite[Proposition 5.2]{AA}; they are faces of $\pp$, hence generalized permutahedra in their own right.

The antipode in $\GP$ was computed by Aguiar and Ardila \cite[Thm.~7.1]{AA} using topological methods:
\begin{equation}\label{eq:antipodegp}
\anti^{\GP}(\pp) = (-1)^{|I|} \sum_{\qq\leq \pp} (-1)^{\codim\qq} \qq.
\end{equation}

The Hopf monoid $\GP_+$ is defined by setting $\GP_+[I]$ to be the $\kk$-vector space spanned by all extended generalized permutahedra in $\Rr^I$.  The product and coproduct are defined in the same way as for $\GP$, so $\GP$ is a Hopf submonoid of $\GP_+$ \cite[\S5.3]{AA}.  

\subsection{The Hopf monoids \texorpdfstring{$\OGP$}{OGP} and \texorpdfstring{$\OGP_+$}{OGP+}}\label{sec:OGP}

\begin{definition}
The \defterm{Hopf monoid of ordered generalized permutahedra} is the Hadamard product $\OGP=\bL^*\x\GP$.  That is, as a vector species, $\OGP[I]=\bL^*[I]\otimes\GP[I]$, and the product and coproduct are defined componentwise on the tensor factors
\end{definition}

Like both $\bL^*$ and $\GP$, the monoid $\OGP$ is commutative but not cocommutative.  We note that $\OGP$ is not linearized as a Hopf monoid (because $\bL^*$ is not), so its antipode is not computable from that of $\GP$ using the methods of Benedetti and Bergeron~\cite{linearized}.

The inclusion $\Mat\inj\GP$ gives rise to an inclusion $\OMat\inj\OGP$, where $\OMat=\bL^*\x\Mat$, the \defterm{Hopf monoid of ordered matroids.}

\begin{theorem}\label{thm:inj}
The symmetrization map $\Symm: \GP\rightarrow \OGP$ defined on $\GP[I]$ by
\[
\Symm(\pp) = \pp^\# = \sum_{w\in\bl[I]} w\otimes\pp
\]
is an injective Hopf morphism. 
\end{theorem}

We also define the \defterm{Hopf monoid of ordered extended generalized permutahedra $\OGP_+$} as follows: $\OGP_+[I]$ is the subspace of $(\bL^*\x\GP_+)[I]$ generated by tensors $\pp\otimes w$ such that $\sigma_w\subset\NN(\pp)$, so that $\pp_w$ is a well-defined vertex of $\pp$.  (The reason for not defining $\OGP_+$ as the full Hadamard product $\bL^*\x\GP_+$ is to avoid basis elements $w\otimes\pp$ such that $\pp$ is unbounded in direction $w$.)

\begin{theorem}
$\OGP_+$ is a submonoid of $\bL^*\x\GP_+$. 
\end{theorem}

\section{Scrope complexes} \label{sec:scrope}

In this section we describe a class of simplicial complexes that will play a key role in the computation of the antipode on $\OGP_+$.

\begin{definition} \label{def:scrope}
A \defterm{Scrope complex} is a simplicial complex on vertices~$[k-1]$ that is either a simplex, or is generated by faces of the form $[k-1]\sm[x,y-1]=[1,x-1]\cup[y,k-1]$, where $1\leq x<y\leq k$.
If $\zz=((x_1,y_1),\dots,(x_r,y_r))$ is a list of ordered pairs of integers in~$[k]$ with $x_i<y_i$ for each $i$, we set $\varphi_i=[k-1]\sm[x_i,y_i-1]$ for $1\leq i\leq r$ and define
\[\Scr(k,\zz) = \left\langle\varphi_1,\dots,\varphi_r\right\rangle.\]
\end{definition}

The facets of a Scrope complex correspond to the intervals $[x,y-1]$ that are minimal with respect to inclusion. By removing redundant generators, we may assume that it is either the full simplex on $[k-1]$, or can be written as $\Scr(k,\zz)$ where $1\leq x_1<\dots<x_r<k$; $1<y_1<\dots<y_r\leq k$; and $x_i<y_i$ for each $i$.

A Scrope complex can be recognized by its facet-vertex incidence matrix, which can be represented as a $r\x(k-1)$ table whose $(i,j)$ entry is $\star$ or $\cdot$ according as $j\in\varphi_i$ or $j\not\in\varphi_i$.  Thus each row consists of a (possibly empty) sequence of dots sandwiched between two (possibly empty) sequences of stars.
For example, if $k=7$ and $\zz=((1,3),(2,4),(3,6),(4,7))$ then $\Scr(n,\zz)=\langle 3456,1456,126,123\rangle$ is represented by the following diagram:
\[\begin{array}{cccccc}
\cdot & \cdot & \star & \star & \star & \star\\
\star & \cdot & \cdot & \star &\star &\star\\
\star & \star & \cdot & \cdot & \cdot & \star\\
\star & \star & \star & \cdot & \cdot & \cdot
\end{array}\]
It is easy to see from this description that the class of Scrope complexes is stable under taking induced subcomplexes.  In fact, an inductive argument shows that Scrope complexes have very simple topology:
\begin{proposition} \label{SG-homotopy}
Every nontrivial Scrope complex is either contractible or a homotopy sphere.
\end{proposition}
\begin{corollary} \label{Euler-SG}
The reduced Euler characteristic of every Scrope complex is 0, 1, or $-1$.
\end{corollary}

\section{The antipode in \texorpdfstring{$\OGP$}{OGP}} \label{sec:antipode}

Our first main result is the following theorem.

\begin{theorem} \label{OGP-antipode}
Let $I$ be a finite set such that $|I|=n$ and $\pp\subseteq\Rr^I$ be a generalized permutahedron of dimension $n-1$,\footnote{We have assumed $\dim\pp=n-1$ for simplicity.  When $\dim\pp<n-1$, the antipode is essentially equivalent but sometimes requires a slight modification; we omit the details.} and let $w$ be a linear ordering of $I$.  Then
\[\anti(w\otimes \pp)
= -\sum_{u\in\bl[I]} (-1)^{\des(u^{-1}w)} u \otimes
\left(\sum_{\substack{\qq\subseteq\pp:\\ \CCC_{\setcomp{w}}\cap\CCC^\circ_Q=\CCC_D}}  \qq
+ \sum_{\substack{\qq\subseteq\pp:\\ D\in\partial\CCC_Q,\ \CCC_{\setcomp{w}}\cap\CCC^\circ_Q\neq\0}} \tilde\chi(\Gamma(Q,w,u)) \qq \right) \label{eq:anti}\]
where $D=D(u,w)$ and $\Gamma(Q,w,u)$ is a certain Scrope complex, described below. This formula is multiplicity-free and cancellation-free.
\end{theorem}

The Scrope complex $\Gamma(Q,w,u)$ is constructed as follows. Let $<$ be the total order on~$I$ defined by $w(1)<w(2)<\cdots<w(n)$.  Let $Q^\natural$ be the naturalization of $Q$ with respect to~$<$, and let $N=N_1|\cdots|N_k$ be the set composition whose blocks are the equivalence classes of~$Q^\natural$. Given a pair $a_i,b_i\in I$ such that $a_i\prec_N b_i$ and $a_i\equiv_D b_i$, let $N_{x_i}$ and $N_{y_i}$ be the blocks of $N$ containing $a_i$ and $b_i$ respectively (so that $x_i<y_i$), and let
\[S_i = N_1 \:|\: \cdots \:|\: N_{x_i-1} \:|\: N_{x_i}\cup\cdots\cup N_{y_i} \:|\: N_{y_i+1} \:|\: \cdots \:|\: N_k.\]
Then the simplicial complex $\Gamma(Q,w,u)$  is a Scrope complex whose vertices correspond to the separators between blocks of $N$.  In the notation of \S\ref{sec:scrope}, we have $\Gamma\isom\Scr(k,\zz)$, where $\zz=((x_1,y_1),\dots,(x_r,y_r))$.  It follows by Proposition \ref{Euler-SG} that $\tilde\chi(\Gamma)\in\{0,-1,1\}$.

\begin{example}
Let $w=u=\mathsf{12345678}$ (as linear orders) so that $D=\hatzero=\mathsf{12345678}$.  Let $Q$ and $Q^\natural$ be the preposets shown below, so that that $N=\mathsf{1|234|567|8}$.
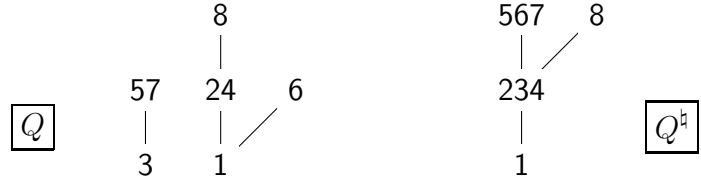
\begin{center}
\begin{tikzpicture}
\newcommand{\sh}{5}
\draw (0,0)--(0,1) (1,2)--(1,0)--(2,1);
							\node[fill=white] at (1,2) {\sf8};
\node[fill=white] at (0,1) {\sf57};		\node[fill=white] at (1,1) {\sf24};		\node[fill=white] at (2,1) {\sf6};
\node[fill=white] at (0,0) {\sf3};		\node[fill=white] at (1,0) {\sf1};
\node at (-1.5,.5) {$\boxed{Q}$};
\begin{scope}[shift={(5,0)}]
\draw (0,0)--(0,2)  (0,1)--(1,2);
\node[fill=white] at (0,2) {\sf567};	\node[fill=white] at (1,2) {\sf8};
\node[fill=white] at (0,1) {\sf234};
\node[fill=white] at (0,0) {\sf1};
\node at (2,.5) {$\boxed{Q^\natural}$};
\end{scope}
\end{tikzpicture}
\end{center}
The simplicial complex $\Gamma=\Gamma(Q,w,u)$ is given by all ordered set partitions that are refined by either $\mathsf{1|234567|8}$ or $\mathsf{1234|567|8}$. The complete list of elements is
\[\{\mathsf{1|234567|8},\ \mathsf{1234|567|8},\ \mathsf{1234567|8},\ \mathsf{1|2345678},\ \mathsf{1234|5678},\ \mathsf{12345678}\}.\]
This is a simplicial complex on the set of separators., consisting of two 1-simplices with a common vertex (so $\tilde\chi(\Gamma)=0$).
\end{example}

\begin{proposition}
Let $\Pi_{n-1}\subset \mathbb{R}^{[n]}$ be the regular permutohedron together with the natural order $w$ on $[n]$. In this case all faces $\qq$ are indexed by ordered set compositions. Applying Theorem \ref{OGP-antipode} we obtain:
\[
\anti(w\otimes \Pi_{n-1})=\sum_{A\textrm{ natural}} \left(\sum_{D(w,u)\refinedbyeq A}u\right)\otimes \qq_A.
\]
\end{proposition}

Locality is illustrated by the fact that only \emph{natural} ordered set compositions appear in this formula.

Meanwhile, since the symmetrization map $\GP\to\OGP$ is a Hopf morphism, the Aguiar--Ardila formula for the antipode in~$\GP$ implies that
\[\sum_{w\in L[I]}\anti^{\OGP}_I(w\otimes \pp) = (-1)^{|I|} \sum_{\substack{\qq\leq\pp\\ u\in L[I]}} (-1)^{\dim \qq}(u\otimes\qq).\]
The right-hand side is cancellation-free and multiplicity-free.  The left-hand side is not cancellation-free, and understanding the cancellation (even for simple examples) is a subtle combinatorial problem.

\section{Hopf monoids of pure ordered complexes}
What conditions on a set of ordered simplicial complexes are required for it to give rise to a Hopf monoid?

\begin{definition} \label{def:HopfClass}
A class $\QQ$ of pure ordered complexes is called a \defterm{Hopf class} if it satisfies the following three conditions.
\begin{enumerate}
\item If $(\Sigma_1,w_1),(\Sigma_2,w_2)\in\QQ$ then $(\Sigma_1*\Sigma_2,w)\in\QQ$ for any $w\in\shuffle(w_1,w_2)$ (where as before $*$ means simplicial join).
\item If $(\Sigma,w)\in\QQ$ and $A$ is an initial segment of $\Sigma$, then $(\Sigma| A,w| A)\in\QQ$. This in particular require that the restriction to all initial segments to be pure.
\item If $(\Sigma,w)\in\QQ$ and $A$ is an initial segment of $\Sigma$, then $(\Sigma/ A,w/ A)\in\QQ$.
\end{enumerate}
\end{definition}

Note that the definition of a Hopf class really requires working with \emph{ordered} complexes --- in the unordered setting, the requirement that all restrictions are pure would limit this definition just to matroid complexes.

\begin{theorem}\label{thm:quasiHopf}
Given a Hopf class $\QQ$, let $\bH_\QQ[I]$ be the linear span of ordered complexes in $\QQ$ on ground set $I$.  Then $\HS[I]$ can be made into a Hopf monoid, with product and coproduct given by
\[\mu_{I,J}(w_1\otimes\Sigma_1)\otimes (w_2\otimes\Sigma_2) = \sum_{w\in \shuffle(w_1, w_2)} w\otimes(\Sigma_1*\Sigma_2),\]
\[\Delta_{I,J}(w\otimes\Sigma) = \begin{cases}
(w|_I\otimes\Sigma|I)\otimes (w/I,\Sigma/I)  & \text{ if 
} w|_I \text{ is an initial segment of } w,\\ 
0 & \text{ otherwise.}\end{cases}.\]
Moreover, the map $\OMat\to\bH_\QQ$ given by $w\otimes M\mapsto w\otimes \II(M)$ is a Hopf monoid monomorphism, where $\II(M)$ is the independence complex of~$M$.
\end{theorem}

\begin{definition}\label{def:decomposable}
An ordered complex $(\Sigma,w)$ on ground set $I$ is \defterm{order-decomposable} if one of the following conditions holds: 
\begin{enumerate}
    \item $\Sigma$ has exactly one facet. 
    \item $\Sigma$ is pure, and for every initial segment with $A$, we have that the complexes $(\Sigma| A, w|A)$ and and $(\Sigma/ A, w/A)$ are pure and order-decomposable.
\end{enumerate}
\end{definition}

The class $\UU$ of all order-decomposable complexes is \textit{universal} in the following sense.

\begin{theorem} \label{universal-Hopf}
The class $\UU$ is a Hopf class. Furthermore, any Hopf class is contained in $\UU$.
\end{theorem}

The following are examples of classes that are contained in $\UU$.

\begin{enumerate}
\item 
The class $\MM$ of ordered matroids is a Hopf class (indeed, it is the prototype for Definition~\ref{def:HopfClass}).

\item
A ordered simplicial complex $(\Sigma,w)$ is \defterm{shifted} if for any face $\sigma\in\Sigma$ and any vertex $e\in\sigma$ then $\sigma\cup{f}\backslash {e}$ is a face for every $f<_we$. That is, replacing any vertex of a face by a smaller one yields another face. The class of pure shifted complexes is not itself a Hopf class since it is not closed under shuffled joins. However, the class of shuffled joins of shifted complexes is a Hopf class. 

\item
Let $(w,M)$ be an ordered matroid.  A \defterm{broken circuit} is obtained by deleting the smallest element of a circuit. The family of all subsets of the ground set that do not contain any broken circuit form a simplicial complex called the \defterm{broken-circuit complex} $BC_w(M)$ (see \cite{bjorner}). These complexes are always pure \cite[Proposition 7.4.2]{bjorner} and lexicographically shellable \cite[7.4.3]{bjorner}.  Broken-circuit complexes are order-decomposable, although the class of broken circuit complexes is not itself a Hopf class (since it is not closed under taking contractions). \\ \\ 
\end{enumerate}

\noindent {\bf Acknowledgements}

We thank Marcelo Aguiar for helping us understand the language of Hopf monoids, and Isabella Novik for supporting a weeklong meeting between the second and third authors in 2015, where this project started. The third author thanks the University of Miami, where much of the work of this project was carried out.
\bibliographystyle{alpha}
\bibliography{biblio}

\end{document}